\documentclass{article}
\usepackage{arxiv}

\usepackage[utf8]{inputenc} 
\usepackage[T1]{fontenc}    
\usepackage{hyperref}       
\usepackage{url}            
\usepackage{booktabs}       
\usepackage{amsfonts}       
\usepackage{nicefrac}       
\usepackage{microtype}      
\usepackage{lipsum}
\usepackage{graphicx}
\usepackage{amsmath}

\title{Coupled Torsional and Transverse Vibration Analysis of Panels Partially Supported by Elastic Beam}

\author{
  Mostafa Bagheri\thanks{mstf.bagheri@gmail.com, http://flyingv.ucsd.edu/mostafa} \\
  Department of Mechanical and Aerospace Engineering\\
  University of California San Diego\\
  La Jolla, CA 15213 \\
  \texttt{mstfbagheri@ucsd.edu} \\
   \And
 Mohammad Mohammadi Aghdam \\
  Department of Mechanical Engineering \\
  Amirkabir University of Technology (Tehran polytechnic) \\
  Tehran 1591634311 \\
  \texttt{aghdam@aut.ac.ir} \\
   \And
 Meitham Amereh \\
  Department of Mechanical Engineering\\
  University of Victoria\\
  Victoria, BC, V8P 5C2 \\
  \texttt{m.amereh@uvic.ca} \\
}

\begin{document}

\maketitle

\begin{abstract}
This study presents torsional and transverse vibration analysis of a solar panel including a rectangular thin plate locally supported by an elastic beam. The plate is totally free in all boundaries, except for the local part attached to the beam. The response of the system, which is subjected to a combination of torsional and transverse vibration, identifies with a couple of PDEs developed by the Euler-Bernoulli assumption and classical plate theory. To calculate the system's natural frequencies, the domain of the solution is discretized by zeroes of the Chebyshev polynomials to apply the Modified Generalized Differential Quadrature method (MGDQ). Furthermore, governing equations along with continuity and boundary conditions are discretized. After obtaining solutions to the eigenvalue problem, several studies are investigated to validate the accuracy of the proposed method. As can be concluded from the tables, MGDQ improves the accuracy of results obtained by GDQ. Results for various case studies reveal that MGDQ is properly devised for the vibration analysis of systems with local boundary and continuity conditions.
\end{abstract}


\section{Introduction}
\label{sec:introduction}
Vibration analysis of plates and panels is an interesting practical subject in both the engineering area and industrial applications. Considerable researches have been performed on the buckling and vibration of plates, however, limited investigations focused more on complicated case studies within the literature. Local suspended plates are among examples of practical subjects in different industries which that need more investigation.

For instance, a suspended free plate on elastic beams can be seen in aerospace and nautical industries. In the area of micro-electromechanical systems (MEMSs), conventional sensors and actuators are among interesting application, where a proof mass attached to one or two suspended beams (e.g. cantilever accelerometers), were studied in \cite{Amereh2016,Tabak2010,Yu2001,Lobontiu2004}. These works only considered free transverse vibration analysis of the system while the torsional part is completely ignored, which plays important role in satellite's solar panels. Solar panels on spacecraft are important samples in which a thin plate is connected to the body of the spacecraft by a thin elastic beam, as shown in Fig. 1.

\begin{figure}[!htbp]
	\centering
	\includegraphics[width=0.65\linewidth]{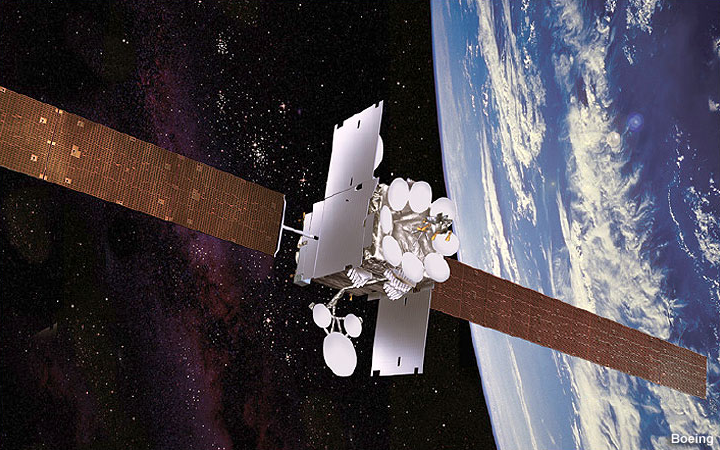}
	\caption{Inmarsat-5 F2 satellite}
	\label{fig:1}
\end{figure}

Notably, that solar panels may experience heavy vibrations due to the environmental effects \cite{Matsuno1996}. Hence, different mode shapes of the vibrating system should be taken into account.

Various techniques are accomplished to find the solution for vibration analysis of shells and plates under specific forms of loading and boundary conditions, series-type method \cite{Rao2007} and integral equations method \cite{Rao1973,Narita1981}, to name a few. Leissa also presented an inclusive and accurate analytical approach for the solution of the free vibration in rectangular plates \cite{Leissa1973,Leissa1984,Rao1973}. Moreover, finite element (FE) and finite strip method (FSM) \cite{Zitnan1996,Chia1985,Leipholz1987}, Rayleigh-Ritz method \cite{Donning1998}, Galerkin method \cite{Belytschko1994}, domain decomposition method \cite{Jang1989,Cheung1989}, energy balance method (EBM) \cite{Sfahani2011}, differential quadrature (DQ) \cite{Bert1993}, and generalized differential quadrature (GDQ) methods \cite{Shu1999,Ng1999} are among different techniques proposed within literature.

The DQ and GDQ methods have been developed by Bellman and Casti \cite{Bellman1971} and Shu and Richards \cite{Shu1990}. The efficiency in computing the weighting coefficients is the principal advantage of GDQ due to its simplicity in the choice of grid points. For the vibration analysis of plates and shells, various studies concerning the validity of DQ/GDQ methods exist in the scientific literature. Most of these studies have investigated typical cases including a combination of standard boundary conditions and uniformly distributed loads. For rectangular plates with different shapes of holes, an independent coordinate coupling method is employed for vibration analysis \cite{Kwak2007}. Additionally, local effects such as the vibration of rectangular plates with concentrated mass–springs \cite{Ingber1992}, distinct elastic edge restraints \cite{Ashour2004}, and thermal nonlinear vibration of the orthotropic DLGS \cite{Arani2013} are investigated within the literature.

This study investigates combined torsional and transverse vibration analysis of a solar panel including a completely free rectangular plate locally supported by an elastic beam. By applying the GDQ method, the domain of the solution together with governing PDEs for the beam and plate are discretized. Based on governing equations along with boundary and continuity conditions, the final eigenvalue problem is developed. The solution to the final equations leads to eigenvalues and eigenvectors which are respectively natural frequencies and mode shapes of the system. To assess the reliability of the proposed method, results for various simplified case studies, in Section 5, are evaluated with available analytical and numerical solutions within the literature. For the complicated cases, results are also validated by commercial FE code, due to lack of information in the literature. Finally, frequency predictions and mode shapes of the general case show an acceptable correlation with finite element results. It is expecting that the results of this paper can be utilized as a benchmark in future studies dealing with the combined transverse and torsional vibration analysis of suspended plates and panels.

\section{Modeling}
\label{sec:odeling}

The system comprises a rectangular isotropic thin elastic plate with dimensions of $a \times b$ and thickness of $h$. The plate is locally attached to a clamped beam with length $l$ and the cross-section of $t\times d$, as shown in Fig. 2. Considering the possibility of resonance in the vibration of the plates, one should pay attention to the material properties and geometry of the system. The plate and the beam are elastic with $E_1$, $E_2$ as modulus of elasticity, and ${\nu}_1$, ${\nu}_2$ as Poisson’s ratio, respectively.
\begin{figure}[!htbp]
	\centering
	\includegraphics[width=0.65\linewidth]{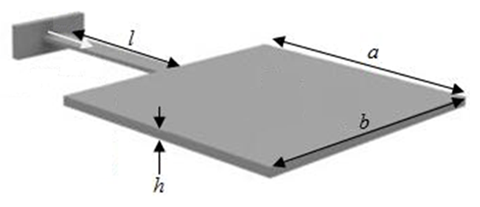}
	\caption{Schematic model}
	\label{fig:2}
\end{figure}

Also, $G_1$ and $I_1$ are shear modulus and polar moment of inertia of the beam, respectively. Due to external excitation, the system may start to vibrate in the form of transverse and/or torsional vibration. Assuming $W(x, y, t)$ as the deflection of the plate, and $U(x^\prime, t)$ and $\theta (x^\prime, t)$ as deflection and rotation of the beam, the governing equations of the system can be obtained as [7]:
\begin{subequations}
	\begin{align} 
	& E_1 I_1 \frac{\partial^4 U(x^\prime,t)}{\partial{x^\prime}^4} + {\rho}_1A_1\frac{\partial^2 U(x^\prime,t)}{\partial t^2}=0 \\
	& G_1 J_1\frac{\partial^2 \theta (x^\prime,t)}{\partial {x^\prime}^2} = \rho_1 I_1\frac{\partial^2\theta (x^\prime,t)}{\partial t^2} \\
	& D\frac{\partial^4 W(x,y,t)}{\partial x^4} + 2D\frac{\partial^4 W(x,y,t)}{\partial x^2\partial y^2} + D\frac{\partial^4 W(x,y,t)}{\partial y^4} + {\rho}_2 h\frac{\partial^2 W(x,y,t)}{\partial t^2} = 0
	\end{align} \label{Eq:eq_1}
\end{subequations}
in which $x^\prime = x + l$ and $A_2 = t \times d$ is the cross-sectional area of the beam.  $\rho_1$ and $\rho_2$ are respectively mass densities of the beam and the plate. Also, $D$ is the flexural rigidity of the plate, which can be determined as
\begin{align} 
    D = \frac{E_1 h^3}{1 2(1 - {\upsilon }^2_2)} \label{Eq:eq_2} 
\end{align}

Assuming $u(x^\prime,t)=U(x^\prime)e^{i\omega^z_bt}$, $\theta (x^\prime,t) = \Theta(x^\prime)e^{i\omega^{\phi}_bt}$, and $w(x\mathrm{,y,}t) = W(x,y)e^{i\omega_pt}$, one may rewrite Eq. \ref{Eq:eq_1} as follows
\begin{subequations}
	\begin{align}
	& \frac{\partial^4 U}{\partial {X^\prime}^4} = {{\overline\omega}^z_b}^2 U \\  
	& \frac{\partial^2 \Theta}{\partial X^2} + {{\overline\omega}^{\phi }_b}^2 \Theta = 0 \\
	& \frac{\partial^4 W}{\partial X^4} + 2\beta^2\frac{\partial^4 W}{\partial X^2 \partial Y^2} + \beta^4\frac{\partial^4 W}{\partial Y^4} = {{\overline\omega}_p}^2 W
	\end{align} \label{Eq:eq_3}
\end{subequations}
where $\omega^z_b$ and ${\omega}^{\phi}_b$ are natural frequencies for transverse and torsional vibration of the beam, respectively, and $\omega_b$ is natural frequency for vibration of plate.

Also, other parameters are defined as
\begin{subequations}
	\begin{align*}
	& {\overline\omega}_b^{z^2} = \frac{\rho_1 A L^4}{E_1 I}\omega_b^{z^2} \\[5pt]
	& {\overline\omega}_b^{\phi^2} = \frac{\rho I L^2}{G J}\omega_b^{\phi^2} \\[5pt]
	& {{\overline\omega}_p}^2 = \frac{\rho_2 h a^4}{D} \omega^2_p \\[5pt]
	& X^\prime = \frac{x^\prime}{L \beta} = \frac{a}{b} \\[5pt]
	& X = \frac{x}{a}, \quad Y = \frac{y}{b}	\end{align*}
\end{subequations}
\setcounter{equation}{3}

Note that the length of the plate is defined as $a$, while $b$ is the plate width. Generally, for a thin plate, shear forces and moments are
\begin{subequations}
	\begin{eqnarray}
	\begin{cases}
	V_x = -D \left( \frac{\partial^3W}{\partial x^3} + (2-v_2)\frac{\partial^3W}{\partial x\partial y^2} \right) \\
	V_y = -D \left( \frac{\partial^3W}{\partial y^3} + (2-v_2)\frac{\partial^3W}{\partial x^2 \partial y} \right)	
	\end{cases} \\
	\left[ \begin{array}{c}
	M_x \\ 
	M_y \\ 
	M_z
	\end{array} \right] = -D \left[ \begin{array}{ccc}
	1 & v_2 & 0 \\ 
	v_2 & 1 & 0 \\ 
	0 & 0 & 1-v_2
	\end{array} \right]
	\left[ \begin{array}{c}
	\frac{\partial^2W}{\partial x^2} \\ 
	\frac{\partial^2W}{\partial y^2} \\ 
	\frac{\partial^2W}{\partial x\partial y}
	\end{array} \right]
	\end{eqnarray} \label{Eq:eq_4}
\end{subequations}

Implementation of boundary conditions, which refers to the free edges of the plate and clamped side of the beam,  as well as continuity of displacements, shear forces, and bending moments at the attached points should be considered. Therefore, boundary equations for the beam section can be written as
\begin{subequations}
	\begin{eqnarray}
	&& u(0,t) = 0 \\[5pt]
	&& \frac{\partial u(0,t)}{\partial x^\prime} = 0 \\[5pt]
	&& \theta (0,t) = 0 
	\end{eqnarray} \label{Eq:eq_5}
\end{subequations}

Also, boundary conditions at free edges of the plate are
\begin{subequations}
	\begin{align}
	& \frac{\partial^2 w(a,y,t)}{\partial x^2} + \nu_2\beta^2\frac{\partial^2 w(a,y,t)}{\partial y^2} = 0 \\[5pt]
	& \frac{\partial^2 w(0,y,t)}{\partial x^2} + \nu_2\beta^2\frac{\partial^2 w(0,y,t)}{\partial y^2} = 0 \quad {\rm for} \ \vert y \vert \ge \frac{d}{2} \\
	& \beta^2 \frac{\partial^2 w(x,0,t)}{\partial y^2} + \nu_2\frac{\partial^2 w(x,0,t)}{\partial x^2} = 0 \\[5pt]
	& \beta^2\frac{\partial^2 w(x,b,t)}{\partial y^2} + \nu_2\frac{\partial^2 w(x,b,t)}{\partial x^2} = 0
	\end{align} \label{Eq:eq_6}
\end{subequations}
and 
\begin{subequations}
	\begin{align}
	& \frac{\partial^3 w(a,y,t)}{\partial x^2}+(2 - \nu_2)\beta^2\frac{\partial^2 w(a,y,t)}{\partial x \partial y^2} = 0 \\[5pt]
	& \frac{\partial^2 w(0,y,t)}{\partial x^2} + (2  - \nu_2)\beta^2\frac{\partial^2 w(0,y,t)}{\partial x \partial y^2} = 0 \quad {\rm for} \ \vert y \vert \ge \frac{d}{2} \\
	& \beta^2\frac{\partial^3 w(x,0,t)}{\partial y^2} + (2 - \nu_2)\frac{\partial^3 w(x,0,t)}{\partial x^2 \partial y} = 0 \\[5pt]
	& \beta^2 \frac{\partial^3 w(x,b,t)}{\partial y^3} + (2 - \nu_2)\frac{\partial^3 w(x,b,t)}{\partial x^2 \partial y} = 0
	\end{align} \label{Eq:eq_7}
\end{subequations}
Continuity conditions at the attached points can be written as
\begin{subequations}
	\begin{align}
	& u(l,t) = w(0,y,t) \\[5pt]  
	& \frac{\partial u(l,t)}{\partial x^\prime} = \frac{\partial w(0,y,t)}{\partial x} \\[5pt] 
	& E_{1} I_{1}\frac{\partial^2 u(l,t)}{\partial {x^\prime}^2} = D \left( \frac{\partial^2 w(0,y,t)} {\partial x^2} + \nu_2 \frac{\partial^2 w(0,y,t)} {\partial y^2} \right) \\  
	& E_1 I_{1}\frac{\partial^3 u(l,t)}{\partial {x^\prime}^2} = D \left( \frac{\partial^2 w(0,y,t)}{\partial x^2} + (2 - \nu_2) \frac{\partial^2 w(0,y,t)}{\partial x \partial y^2} \right)
	\end{align} \label{Eq:eq_8}
\end{subequations}

in which the domain of $y$ in Eqs. \ref{Eq:eq_8} is $\vert y \vert \ge \frac{d}{2}$.

\section{GDQ Implementation}
To commence the procedure of GDQ method, the domain of the solution should be first discretized into several grid points. Zeros of the Chebyshev polynomials is one of the best choices to discretize the domain, see \cite{Shu2001}. As a result, the plate is discretized into $N \times M$ grid points in $x$ and $y$ directions, respectively, while the beam may be divided into three or five rows of $S$ grid points through $x$ direction, as shown in Fig. 3. Afterward, governing partial differential equations are expanded over these grid points. By assuming three or five rows of grid points in the width of the beam, continuity conditions can be satisfied.
\begin{figure}[!htbp]
	\centering
	\includegraphics[width=0.65\linewidth]{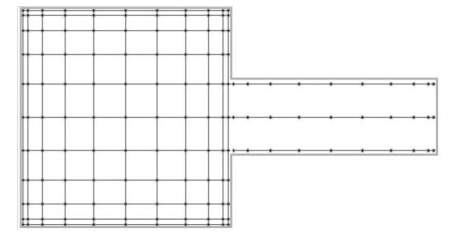}
	\caption{Discretization of the system into grid points}
	\label{fig:3}
\end{figure}

According to the GDQ method, a general form of the approximation for the calculation of higher-order derivatives is
\begin{eqnarray} 
u^m_x (x_i,t) = \sum^b_{j = 1} \tau^{(m)}_{i j} u(x_j,t) \qquad {\rm for} \quad i = 1,\cdots, S
\label{Eq:eq_9}
\end{eqnarray}
where $u^m_x (x_i,t)$ is the  $m^{\rm th}$ derivative of beam deflection at point $(x_i,b)$ is the number of grid points and $\tau^{(m)}_{ij}$ is the weighting coefficients for the $m^{\rm th}$ derivative calculation. The grid points on the beam and plate based on zeros of the Chebyshev polynomials can be determined as:
\begin{subequations}
	\begin{align} 
	x^\prime_k & = \frac{L}{2} \left[1-{{\rm cos} \left(\frac{k-1}{S-1} \pi \right)}\right] \\[5pt]
	x_i & = \frac{a}{2} \left[1-{{\rm cos} \left(\frac{i-1}{N-1} \pi \right)\ }\right] \\[5pt]
	y_j &= \frac{b}{2} \left[1-{{\rm cos} \left(\frac{j-1}{M-1} \pi \right)\ }\right]
	\end{align} \label{Eq:eq_10}
\end{subequations}

In discretization, Eq. \eqref{Eq:eq_11} gives the weighting coefficients for higher order derivatives
\begin{subequations}
	\begin{align} 
	& \zeta^{(m)}_{rk} = m \left( \eta^{(m-1)}_{rr}\eta^{(1)}_{rk} - \frac{\eta^{(m - 1)}_{rk}}{{x^\prime}_r - {x^\prime}_k} \right) \quad \qquad {\rm for} \quad r \neq k \quad r, k = 1, \cdots, S \\
	& c^{(m)}_{pi} = m \left( c^{(m-1)}_{pp}c^{(1)}_{pi}-\frac{c^{(m-1)}_{pi}}{x_p-x_i} \right) \quad \qquad {\rm for} \quad r \neq k \quad r, k = 1, \cdots, S \\
	& {\overline{c}}^{(m)}_{lj} = m \left( {\overline{c}}^{(m - 1)}_{ll}{\overline{c}}^{(1)}_{lj} - \frac{{\overline{c}}^{(m-1)}_{lj}}{y_l - y_j} \right) \quad \qquad {\rm for} \quad r \neq k \quad r, k = 1, \cdots, S
	\end{align} \label{Eq:eq_11}
\end{subequations}

The $\zeta^{(m)}_{rk}$, $c^{(m)}_{pi}$, and ${\overline{c}}^{(m)}_{lj}$ can be respectively derived based on the $(m-1)^{\rm th}$ order coefficients $\zeta^{(m-1)}_{rk}$, $c^{(m-1)}_{pi}$, and ${\overline{c}}^{(m-1)}_{lj}$. The $\zeta^{(m)}_{rk}$, $c^{(m)}_{pi}$, and ${\overline{c}}^{(m)}_{lj}$ can be obtained from:

\begin{subequations}
	\begin{align} 
	\zeta^{(m)}_{rr} & = \sum^S_{\begin{array}{c} k = 1 \\ k \neq r \end{array}} \zeta^{(m)}_{rk} \qquad {\rm for \ } \  r = 1, \cdots, S \\[5pt]
	c^{(m)}_{pp} & = \sum^N_{\begin{array}{c} l = 1 \\ l \neq p \end{array}} c^{(m)}_{pl} \qquad {\rm for \ } \  p = 1, \cdots, N \\[5pt]
	{\overline{c}}^{(m)}_{ll} & = \sum^M_{\begin{array}{c} j = 1 \\ j \neq l \end{array}} {\overline{c}}^{(m)}_{ij} \qquad {\rm for \ } \ l = 1, \cdots, M
	\end{align} \label{Eq:eq_12}
\end{subequations}

Based on the recursive formulas \ref{Eq:eq_11} and \ref{Eq:eq_12}, ${\eta }^{(1)}_{rk}$, $c^{(1)}_{pi}$ and ${\overline{c}}^{(1)}_{lj}$ can be calculated as follow:
\begin{subequations}
	\begin{align} 
	& \eta^{(1)}_{rk} = \frac{\mathrm{P}({x^\prime}_r)}{({x^\prime}_r - {x^\prime}_k)\mathrm{P}({x^\prime}_k)} \qquad {\rm for \ } \ r \neq k \\[5pt]  
    & c^{(1)}_{pi} = \frac{Q(x_p)}{(x_p-x_i)Q(x_i)} \qquad {\rm for \ } \ p \neq i \\[5pt] 
	& {\overline{c}}^{(1)}_{lj} = \frac{H(y_l)}{(y_l - y_j)H(y_j)} \qquad {\rm for \ } \ l \neq j 
	\end{align} \label{Eq:eq_13}
\end{subequations}
where
\begin{align}
    & P(x^\prime_r) = \prod^S_{ \begin{array}{c}
    	k = 1 \\ 
    	k \neq r \end{array}
    }{ \left( {x^\prime}_r - {x^\prime}_k \right) }, \qquad Q(x_p) = \prod^N_{ \begin{array}{c}
    i = 1 \\ 
    i\neq p \end{array}
    }{ \left( x_p - x_i \right) }, \qquad H(x_l) = \prod^M_{ \begin{array}{c}
    j = 1 \\ 
    j \neq l \end{array}}{ \left( x_l - x_j \right) } \nonumber
\end{align}

Also, one may obtained $\eta^{(1)}_{kk}$, $c^{(1)}_{ii}$ and ${\overline{c}}^{(1)}_{jj}$ using:
\begin{subequations}
	\begin{align} 
	{\eta }^{(1)}_{kk} & = - \sum^S_{ \begin{array}{c}
		k = 1 \\ 
		k \neq r \end{array}
	}{{\eta }^{(1)}_{rk}} \qquad {\rm for \ } \ r = 1, 2,\cdots, N \\ 
	c^{(1)}_{ii} & = -\sum^N_{ \begin{array}{c}
		i = 1 \\ 
		i \neq p \end{array}
	}{c^{(1)}_{pi}} \qquad {\rm for \ } \ p = 1, 2, \cdots, N \\  
	{\overline{c}}^{(1)}_{jj} & = - \sum^M_{ \begin{array}{c}
		j = 1 \\ 
		j \neq l \end{array}
	}{{\overline{c}}^{(1)}_{lj}} \qquad {\rm for \ } \ l = 1, 2, \cdots, N
	\end{align} \label{Eq:eq_14}
\end{subequations}

According to the definition of the derivative terms and applying the GDQ approximation, Eq. \eqref{Eq:eq_3} can be rewritten as
\begin{subequations}
	\begin{align}
	& \sum^S_{k = 1}{\zeta^{(4)}_{rk} U_k} = {\overline\omega}^2 U_r \\ 
	& \sum^S_{k = 1}{\zeta^{(2)}_{rk} \Theta_k} = - {{\overline\omega}^{\phi}_b}^2 \Theta_r \quad {\rm for} r = 3,\cdots,S-2 \\
	& \sum^{N}_{k = 1}{c^{(4)}_{ik}W_{kj}} + 2 \beta^2 \sum^{M}_{m = 1}{{\overline{c}}^{(2)}_{jm}}\sum^{N}_{k = 1}{c^{(2)}_{ik} W_{km}} + \beta^4 \sum^{M}_{m = 1}{{\overline{c}}^{(4)}_{jm} W_{im}} = {{\overline\omega}_p}^2 W_{ij}
	\end{align} \label{Eq:eq_15}
\end{subequations}

for $i = 3, \cdots, N-2$ and $j = 3, \cdots, M-2$. It is worth mentioning that the aforementioned equations should be considered for all grid points except boundary and adjacent nodes, see \cite{Shu1997}. Boundary conditions for the beam section at the fixed end, which have been written in Eq. \eqref{Eq:eq_5}, can be respectively discretized as $U_1 = 0$, $\sum^S_{k = 1}\zeta^{(1)}_{1k} U_k = 0$, and $\Theta_1 = 0$.

In addition to that, the discretized forms of the boundary conditions of zero normal moment at free edges in Eqs. \eqref{Eq:eq_6} and \eqref{Eq:eq_7} are

\begin{subequations}
	\begin{align}
	& \sum^{N}_{k = 1}{c^{(2)}_{Nk} W_{kj}} + \nu_2 \beta^2 \sum^{M}_{m = 1} {{\overline{c}}^{(2)}_{jm} W_{Nm}} = 0 \\
	& \sum^{N}_{k = 1}{c^{(2)}_{1k}W_{kj}} + \nu_2 \beta^2 \sum^{M}_{m = 1} {{\overline{c}}^{(2)}_{jm} W_{1m}} = 0 \qquad {\rm for \ } j = 3, 4, \cdots, M-3, M-2, {\rm and} \nonumber \\
	& \beta^2 \sum^M_{m = 1} {{\overline{c}}^{(2)}_{1m} W_{im}} + \nu_2 \sum^{N}_{k = 1}{c^{(2)}_{ik} W_{k1}} = 0 \\
	& \beta^2 \sum^{M}_{m = 1} {{\overline{c}}^{(2)}_{Mm} W_{im}} + \nu_2 \sum^{N}_{k = 1}{c^{(2)}_{ik} W_{kM}} = 0 \qquad {\rm for \ } i = 1, 2, \cdots ,N \nonumber
	\end{align} \label{Eq:eq_16}
\end{subequations}

Also, we have
\begin{subequations}
	\begin{align}
	& \sum^N_{k = 1}{c^{(3)}_{Nk}W_{kj}} + (2 - \nu_2){\beta }^2\sum^{N}_{k = 1}{c^{(1)}_{Nk}}\sum^{M}_{m = 1}{{\overline{c}}^{(2)}_{jl}W_{km}} = 0 \\ 
	& \sum^N_{k = 1}{c^{(3)}_{1k}W_{kj}} + (2 - \nu_2) \beta^2\sum^N_{k = 1} {c^{(1)}_{1k}} \sum^{M}_{m = 1} {{\overline{c}}^{(2)}_{jl}W_{km}} = 0
	\qquad {\rm for \ } j = 3, 4, \cdots, M-3, M-2, {\rm and} \nonumber \\
	& \beta^2\sum^M_{m = 1} {{\overline{c}}^{(3)}_{1m} W_{im}} + (2 - \nu_2) \sum^N_{k = 1} {c^{(2)}_{1m}}\sum^M_{m = 1} {{\overline{c}}^{(1)}_{ik} W_{km}} = 0 \\ 
	& \beta^2\sum^M_{m = 1} {{\overline{c}}^{(3)}_{Mm}W_{im}} + (2 - \nu_2) \sum^N_{k = 1} c^{(2)}_{Mm} \sum^M_{m = 1} {{\overline{c}}^{(1)}_{ik}W_{km}} = 0 \qquad {\rm for \ }  i = 1, 2, \cdots, N \nonumber
	\end{align} \label{Eq:eq_17}
\end{subequations}

Moreover, the continuity condition for the local part (Eq. \eqref{Eq:eq_8}) can be written as follows
\begin{subequations}
	\begin{align}
	& U_S = W_{1j} \\ 
	& \sum^S_{r = 1} {\zeta^{(1)}_{Sr} U_i} = \sum^N_{k = 1} {c^{(1)}_{1k} W_{kj}} \\  
	& E_{1}I_{1} \sum^S_{r = 1} {\zeta^{(2)}_{Sr}U_r} = -D \left( \sum^N_{k = 1} {c^{(2)}_{1k}W_{\mathrm{kj}}} + \nu_2 \sum^M_{m = 1} {{\overline{c}}^{(2)}_{jm} W_{1m}} \right) \\
	& E_{1}I_{1}\sum^S_{r = 1}{\zeta^{(3)}_{Sr}U_r} = -D \left( \sum^N_{k=1}{c^{(2)}_{1k} W_{kj}} + (2 - \nu_2) \sum^N_{k = 1}{c^{(1)}_{1k}} \sum^M_{m = 1} {{\overline{c}}^{(2)}_{jm} W_{km}} \right) 
	\end{align} \label{Eq:eq_18}
\end{subequations}
in which the domain of $k$ in Eqs. \ref{Eq:eq_18} is $(\frac{b-d}{2})(\frac{N-1}{b})+1 \le j \le (\frac{b+d}{2})(\frac{N-1}{b})+1$.

\section{Solution procedure} \label{Solution procedure}
Considering governing Eqs. \eqref{Eq:eq_15}, boundary conditions \eqref{Eq:eq_16} and \eqref{Eq:eq_17} and continuity conditions \eqref{Eq:eq_18}, one can write these systems of equations in form of two sets of algebraic equations as following:
\begin{subequations}
	\begin{align}
	& \left[ A_{IB} \right] \left\{ W_B \right\} + \left[ A_{II} \right] \left\{ W_I \right\} = {\Omega}^2 \left\{ W_I \right\} \label{Eq:eq_19a} \\[5pt]
	& \left[ A_{BB} \right] \left\{ W_B \right\} + \left[ A_{BI} \right] \left\{ W_I \right\} = 0 \label{Eq:eq_19b} 
	\end{align} \label{Eq:eq_19}
\end{subequations}
where $\{W_B\}$ and $\{W_I\}$ are nodes deflection inside the domain and on boundaries, respectively. Furthermore, dimensions of $[A_{I I}]_{P \times P}$, $[A_{B B}]_{Q \times Q}$, $[A_{I B}]_{P \times Q}$ and $[A_{BI}]_{Q \times P}$ are defined such that $P = (S - 4) + (M - 4) \times (N - 4)$ and $Q = (4M + 4N - 12)$.

Substitution of Eq. \eqref{Eq:eq_19a} in \eqref{Eq:eq_19b} gives the final system of eigenvalue equation as:
\begin{align} 
    \left\{ \left[ A_{II} \right] - \left[ A_{IB} \right] {[A_{BB}]}^{-1} \left[ A_{BI} \right] \right\} \{W_I\} = {\omega}^2 \{W_I\} \label{Eq:eq_20}
\end{align}

The standard form of eigenvalue equation can be obtained by restating Eq. \eqref{Eq:eq_20} as
\begin{eqnarray} 
\left[ \left\{ \left[ A_{II} \right] - \left[ A_{IB} \right] {[A_{BB}]}^{-1} \left[ A_{BI} \right] \right\} - I \right] \{W_I\} = {\omega}^2 \left\{ W_I t \right\} \label{Eq:eq_21}
\end{eqnarray}
where $\{I\}$ is the identity matrix. After implementing the values in stiffness matrix in Eq. (21), by solving this standard eigenvalue equation, natural frequencies of the system together with related mode shapes could be determined.

\section{Results and discussions}
Free torsional and transverse vibration analysis of a system comprises a free plate locally supported by an elastic beam is presented. At first, three types of equations, i.e. governing equations, equations of boundary conditions and equations of continuity conditions are developed. To solve these equations, GDQ method is applied to the domain of the solution. Accordingly, both domains of the solution and governing equations, together with boundary/continuity equations are discretized. Thereafter, all equations are rewritten in the form of two sets of equations giving the final system of eigenvalue equation. The system is eventually simplified into several widely-known cases of vibration analysis, e.g. completely free (FFFF) plate, cantilever (CFFF) plate, clamped-free beam and beam with concentrated mass, to determine the accuracy of the proposed solution, the system is simplified into several. Also, the general model has been investigated in the final case in which the accuracy of the results is evaluated. Table 1 shows the properties of the plate and beam.

\begin{table}[htbp]
	\caption{Material and geometry properties}
	\vspace{0.2cm}
	\centering
	\begin{tabular}{|c|c|} \hline 
		\textbf{Properties} & \textbf{Values} \\ \hline 
		$a$ & $1 \ m$ \\ \hline 
		$b$ & $1 \ m$ \\ \hline 
		$L$ & $1 \ m$ \\ \hline 
		$d$ & $0.1 \ m$ \\ \hline 
		$t$ & $0.5 \ cm$ \\ \hline
		$\nu_1, \nu_2$ & $3 \ cm$ \\ \hline 
		$\rho_1, \rho_2$ & $2330\ kg/m^3$ \\ \hline 
		$E_1, E_2$ & $200 \ GPa$ \\ \hline
	\end{tabular}
\end{table}

\textbf{Case study 1}
As the first case study, material properties and geometry parameters for the plate are assumed to be the same as those of the beam which simplifies the model to a simple cantilever beam. Hence, identical width, thickness and material properties should be considered, i.e. $t=h, d=2b, E_1 = E_2, \rho_1 = \rho_2$. Moreover, the Poisson ratio of the plate section should be neglected to eliminate the deformation of the plate in the $y$ direction. Taking into account the aforementioned properties, the model is simplified to a homogeneous clamped--free beam with length$\ L^*=L+a$. Due to the attachment of the plate to the end of the beam, vibration of the plate may cause torsional vibration over the beam. Therefore, 2D discretization is applied through the domain of the beam. It is worth mentioning that analytical solutions can be found for natural frequencies of both transverse and torsional vibrations of clamped-free beams using modal analysis as \cite{Rao2007}:
\begin{subequations}
	\begin{align}
	& {\rm cos} \beta L^* \times \mathrm{cosh} \beta L^* + 1 = 0 \\[5pt] 
	& \mathrm{cos} \frac{\omega^{\theta}_j L^*}{c} = 0  
	\end{align} \label{Eq:eq_22}
\end{subequations} 
in which
\begin{subequations}
    \begin{align}
    & \beta^4 = \frac{\rho A({\omega^z_i)}^2}{EI}, \qquad c = \sqrt{\frac{G}\rho}, \qquad \omega_j^\theta = \frac{\alpha_j c}{l} \nonumber
    \end{align}
\end{subequations}

Note that $\omega^z_i$ and $\omega^{\theta }_j$ are respectively the $i^{\rm th}$ natural frequency of transverse vibration and $j^{\rm th}$ natural frequency of torsional vibration of the beam.

Tables 2 and 3 respectively demonstrate the first eight natural frequencies of the system for both transverse and torsional vibration, computed by both GDQ and frequency Eqs. \eqref{Eq:eq_22}.
\begin{table*}[bp]
	\caption{First eight natural frequencies of transverse vibration of clamped-free beam}
	\vspace{0.2cm}
	\centering
	\begin{tabular}{|c|c|c|c|c|c|c|c|c|} \hline 
		\textbf{Method/Grid Size} & $\Omega_1$ & $\Omega_2$ & $\Omega_3$ & $\Omega_4$ & $\Omega_5$ & $\Omega_6$ & $\Omega_7$ & $\Omega_8$ \\ \hline 
		\textbf{Exact Method (Eq. 22.a)} & 1.875 & 4.694 & 7.855 & 10.996 & 14.137 & 17.279 & 20.420 & 23.562 \\ \hline 
		\textbf{$5 \times 9$} & 1.845 & 4.569 & 7.736 & 10.825 & 13.984 & 17.259 & 20.390 & 23.493 \\ \hline 
		\textbf{$5 \times 11$} & 1.863 & 4.679 & 7.839 & 10.971 & 14.120 & 17.264 & 20.402 & 23.540 \\ \hline 
		\textbf{$5 \times 15$} & 1.875 & 4.695 & 7.843 & 10.982 & 14.129 & 17.268 & 20.407 & 23.549 \\ \hline 
		\textbf{Error (\%)} & 0.0320 & 0.0128 & 0.1515 & 0.1255 & 0.0559 & 0.0602 & 0.0686 & 0.0531 \\ \hline 
	\end{tabular}
\end{table*}

Another simple case can be considered by assuming a plate with large density and small sizes which reduces the model to a simple cantilever beam with concentrated end mass. To achieve this, the plate is assumed to be made of Lead with density equals to $\rho_1 = 11.34\ gr/cm^3$ and geometry properties as $a = b = h = 0.1m$. Again, the frequency equation for transverse and torsional vibration of this simple case can be found analytically as \cite{Rao2007}:
\begin{subequations}
	\begin{align}
	& 1 + \frac{1}{ {\rm cos} (\beta l) \ {\rm cosh} (\beta l) } - R_z \  \beta l \left( {{\rm tan} (\beta l) - {\rm tanh} (\beta l)} \right) = 0 \\[5pt]
	& \beta^4 = \frac{\rho A({\omega^z_i)}^2}{E I} \\[3pt] 
	& \alpha_j \ {\rm tan} (\alpha_j) = R_\theta \\[5pt]
	& \omega^{\theta}_j = \frac{\alpha_j c}{l}  
	\end{align} \label{Eq:eq_23}
\end{subequations} 

\begin{table*}[htbp]
	\caption{First eight natural frequencies of torsional vibration of clamped-free beam }
	\vspace{0.2cm}
	\centering
	\begin{tabular}{|c|c|c|c|c|c|c|c|c|} \hline 
		\textbf{Method/Grid Size} & $\Omega_1$ & $\Omega_2$ & $\Omega_3$ & $\Omega_4$ & $\Omega_5$ & $\Omega_6$ & $\Omega_7$ & $\Omega_8$ \\ \hline 
		\textbf{Exact Method (Eq. 22.b)} & 1.571 & 4.712 & 7.854 & 10.996 & 14.137 & 17.279 & 20.420 & 23.562 \\ \hline 
		\textbf{$5 \times 9$} & 1.571 & 4.715 & 7.875 & 11.096 & 14.178 & 17.2910 & 20.462 & 23.583 \\ \hline 
		\textbf{$5 \times 11$} & 1.571 & 4.714 & 7.864 & 10.796 & 14.158 & 17.289 & 20.442 & 23.573 \\ \hline 
		\textbf{$5 \times 15$} & 1.571 & 4.710 & 7.852 & 10.988 & 14.128 & 17.281 & 20.418 & 23.573 \\ \hline 
		\textbf{Error (\%)} & 0.013 & 0.045 & 0.024 & 0.071 & 0.067 & 0.014 & 0.011 & 0.045 \\ \hline 
	\end{tabular}
\end{table*}
where $R_z = \frac{M}{\rho A l}$  is the ratio of the attached mass $M$ to the mass of the beam $\rho A l$. Also, $R_\theta$ is equal to $\frac{\rho j l}{I_d}$ in which $I_d$ is the moment of inertia of the plate and $\rho, j$ and $l$ are properties of the beam. Based on the proposed geometry and material properties, one may obtain $R_z = 9.734$ and $R_\theta = 0.051$. Similarly, Tables 4 and 5 contain first eight natural frequencies of the beam with concentrated mass, obtained by GDQ and Eqs. \ref{Eq:eq_23}. Results show reasonably accurate predictions for both cases.
\begin{table*}[bp]
	\caption{First eight natural frequencies of clamped-free beam with concentrated mass $R_z = 9.734$; transverse vibration }
	\vspace{0.2cm}
	\centering
	\begin{tabular}{|c|c|c|c|c|c|c|c|c|} \hline 
		\textbf{Method / Grid Size} & \multicolumn{8}{|p{3.7in}|}{\bf Transverse Frequencies: $R_z = 9.734$} \\ \hline 
		& $\Omega_1$ & $\Omega_2$ & $\Omega_3$ & $\Omega_4$ & $\Omega_5$ & $\Omega_6$ & $\Omega_7$ & $\Omega_8$ \\ \hline 
		\textbf{Exact Method (Eq. 23.a)} & 0.741 & 1.571 & 3.939 & 4.712 & 7.076 & 7.854 & 10.215 & 10.996 \\ \hline 
		\textbf{$5 \times 9$} & 0.763 & 1.597 & 3.948 & 4.733 & 7.284 & 7.963 & 10.321 & 10.962 \\ \hline 
		\textbf{$5 \times 11$} & 0.726 & 1.586 & 3.929 & 4.724 & 7.176 & 7.859 & 10.319 & 10.977 \\ \hline 
		\textbf{$5 \times 15$} & 0.740 & 1.571 & 3.938 & 4.709 & 7.069 & 7.853 & 10.208 & 10.989 \\ \hline 
		\textbf{Error (\%)} & 0.0675 & 0.0255 & 0.0102 & 0.0785 & 0.0862 & 0.0115 & 0.0705 & 0.0591 \\ \hline 
	\end{tabular}
\end{table*}

\begin{table*}[htbp]
	\caption{First eight natural frequencies of clamped-free beam with concentrated mass $R_\theta = 0.051$; torsional vibration}
	\vspace{0.2cm}
	\centering
	\begin{tabular}{|c|c|c|c|c|c|c|c|c|} \hline 
		\textbf{Method / Grid Size} & \multicolumn{8}{|p{3.6in}|}{\bf Torsional Frequencies; $R_\theta = 0.051$} \\ \hline 
		& $\Omega_1$ & $\Omega_2$ & $\Omega_3$ & $\Omega_4$ & $\Omega_5$ & $\Omega_6$ & $\Omega_7$ & $\Omega_8$ \\ \hline 
		\textbf{Exact Method (Eq. 23.b)} & 1.571 & 3.158 & 4.712 & 6.291 & 7.854 & 9.430 & 10.996 & 12.571 \\ \hline 
		\textbf{$5 \times 9$}\textit{} & 1.573 & 3.160 & 4.715 & 6.294 & 7.869 & 9.459 & 11.054 & 12.549 \\ \hline 
		\textbf{$5 \times 11$} & 1.572 & 3.159 & 4.714 & 6.293 & 7.867 & 9.449 & 11.005 & 12.560 \\ \hline 
		\textbf{$5 \times 15$} & 1.570 & 3.157 & 4.710 & 6.290 & 7.849 & 9.421 & 10.985 & 12.562 \\ \hline 
		\textbf{Error (\%)} & 0.0382 & 0.0443 & 0.0552 & 0.0286 & 0.0701 & 0.1007 & 0.0946 & 0.0684 \\ \hline 
	\end{tabular}
\end{table*}

\textbf{Case study 2}
As another simple case for evaluation of the model and solution procedure, one may ignore the effects of the beam by assuming zero geometry and material parameters which reduces the system to a completely free plate. As a result, the beam exerts no force and moment on the plate. Table 6 includes the results of first five natural frequencies obtained by GDQ method. Also, table 6 shows the results of direct GDQ method obtained by Shu and Du \cite{Shu1997}, analytical results presented by Leissa \cite{Leissa1973,Leissa1984} and also results obtained using commercial FE code ANSYS. Predictions show good correlations with other methods. It is worth mentioning that by reselecting grid points by modified formulation provided in \cite{Shu1997}, the accuracy of the results is significantly improved.
\begin{table*}[htbp]
	\caption{First five natural frequencies of FFFF plate}
	\vspace{0.2cm}
	\centering
	\begin{tabular}{|c|c|c|c|c|c|c|} \hline 
		\textbf{Method / Grid Size} & $\Omega_1$ & $\Omega_2$ & $\Omega_3$ & $\Omega_4$ & $\Omega_5$ & $\Omega_6$ \\ \hline 
		\textbf{Leissa \cite{Leissa1973}} & 13.489 & 19.789 & 24.432 & 35.024 & 35.024 & 61.526 \\ \hline 
		\textbf{Leissa \& Narita \cite{Leissa1984}} & 13.468 & 19.596 & 24.271 & 34.801 & 34.801 & 61.111 \\ \hline 
		\textbf{Shu and Du (15*15) \cite{Shu1997}} & 10.303 & 19.596 & 22.146 & 30.026 & 30.803 & - \\ \hline 
		\textbf{Shu and Du${}^{ }$(12*12) \cite{Shu1997}} & 13.454 & 19.597 & 24.271 & 34.815 & 34.817 & - \\ \hline 
		\textbf{FEM} & 13.461 & 19.665 & 24.289 & 34.912 & 34.825 &  \\ \hline 
		\textbf{GDQ} &  &  &  &  &  &  \\ \hline 
		\textbf{$9 \times 9$} & 10.934 & 19.365 & 21.935 & 29.432 & 29.536 & - \\ \hline 
		\textbf{$11 \times 11$} & 10.639 & 19.685 & 22.854 & 31.342 & 30.762 & - \\ \hline 
		\textbf{$15 \times 15$} & 10.303 & 19.596 & 22.146 & 30.026 & 30.803 & - \\ \hline 
		\textbf{Error (\%)} & 23.619 & 0.975 & 9.356 & 14.270 & 12.051 &  \\ \hline 
		\textbf{Modified GDQ} &  &  &  &  &  &  \\ \hline 
		\textbf{$9 \times 9$} & 14.065 & 19.968 & 24.696 & 34.938 & 34.645 & - \\ \hline 
		\textbf{$11 \times 11$} & 13.164 & 19.492 & 24.434 & 34.884 & 34.783 & - \\ \hline 
		\textbf{$15 \times 15$} & 13.475 & 19.598 & 24.268 & 34.832 & 34.828 & - \\ \hline 
		\textbf{Error (\%)} & 0.103 & 0.965 & 0.671 & 0.548 & 0.559 &  \\ \hline 
	\end{tabular}
\end{table*}

\textbf{Case study 3}
As the final test case to evaluate the performance of the presented method, the beam is considered to be rigid with identical widths for both the beam and plate. This will reduce the system to a cantilever plate. To this end, the Young and shear modulus of elasticity of the beam should tend to infinity. Table 7 includes the GDQ predictions for the first five natural frequencies of the plate. Also, the result of direct GDQ and the analytical method by Leissa \cite{Leissa1973,Leissa1984} are presented in Table 7. To prove the validity of the presented method, results of FE code (ANSYS) are also included in the table. To develop the accuracy of the method, the same procedure as the previous case study may be applied. The new results are presented in Table 7 as well.
\begin{table*}[tp]
	\caption{First five natural frequencies of CFFF plate}
	\vspace{0.2cm}
	\centering
	\begin{tabular}{|c|c|c|c|c|c|c|} \hline 
		\textbf{Method / Grid Size} & $\Omega_1$ & $\Omega_2$ & $\Omega_3$ & $\Omega_4$ & $\Omega_5$ & $\Omega_6$ \\ \hline 
		\textbf{Leissa \cite{Lobontiu2004}} & 3.492 & 8.525 & 21.429 & 27.331 & 31.111 & 54.443 \\ \hline 
		\textbf{Shu and Du (15*15) \cite{Shu1997}} & 3.898 & 9.459 & 20.206 & 26.150 & 26.500 & - \\ \hline 
		\textbf{Shu and Du (12*12) \cite{Shu1997}} & 3.485 & 8.604 & 21.586 & 27.230 & 31.358 & - \\ \hline 
		\textbf{FEM} & 3.481 & 8.502 & 21.456 & 27.401 & 30.021 &  \\ \hline 
		\textbf{GDQ} &  &  &  &  &  &  \\ \hline 
		\textbf{$9 \times 9$} & 3.486 & 8.616 & 19.894 & 25.965 & 27.925 & - \\ \hline 
		\textbf{$11 \times 11$} & 3.768 & 8.914 & 20.065 & 26.425 & 29.682 & - \\ \hline 
		\textbf{$15 \times 15$} & 3.898 & 9.459 & 20.206 & 26.150 & 26.500 & - \\ \hline 
		\textbf{Error (\%)} & 11.626 & 10.956 & 5.707 & 4.321 & 14.821 &  \\ \hline 
		\textbf{Modified GDQ} &  &  &  &  &  &  \\ \hline 
		\textbf{$9 \times 9$} & 3.467 & 8.725 & 21.104 & 26.935 & 31.825 & - \\ \hline 
		\textbf{$11 \times 11$} & 3.482 & 8.625 & 21.238 & 27.189 & 31.539 & - \\ \hline 
		\textbf{$15 \times 15$} & 3.495 & 8.564 & 21.462 & 27.312 & 31.261 & - \\ \hline 
		\textbf{Error (\%)} & 0.085 & 0.457 & 0.154 & 0.069 & 0.482 &  \\ \hline 
	\end{tabular}
\end{table*}

\textbf{Case study 4:}
The final case is the complicated system consists of a regionally suspended plate connected to the elastic beam. The standard form of eigenvalue equations for the combined transverse/torsional vibration of the system is solved and consequently eigenvalues and eigenvectors are derived. To validate the calculation procedure, finite element code ANSYS was also used to analyze the system. Figures 4 to 8 depicts five different mode shapes of the system together with associated frequencies, which are predicted by the presented GDQ technique. As shown in these figures, the system can vibrate based on various mode shapes of the beam and the plate and also their combination. For instance, Figs. 4 and 5 depict mode shapes of the system once the first modes of transverse and torsional vibrations of the beam are excited, respectively. Besides, Figs. 6 to 8 show mode shapes of the system as the modes of beam and plate are simultaneously excited. Included in these figures are also results of FE analysis for both frequencies and mode shapes of the system, which indicates reasonable accuracy for the proposed approach. The range of errors in predictions is about $4.4\%$ to $7.8\%$, in which minimum and maximum errors belong to predicted frequencies of $2.675$ Hz and $147.494$ Hz, respectively.

\begin{figure}[!htbp]
	\centering
	\includegraphics[width=0.5\linewidth]{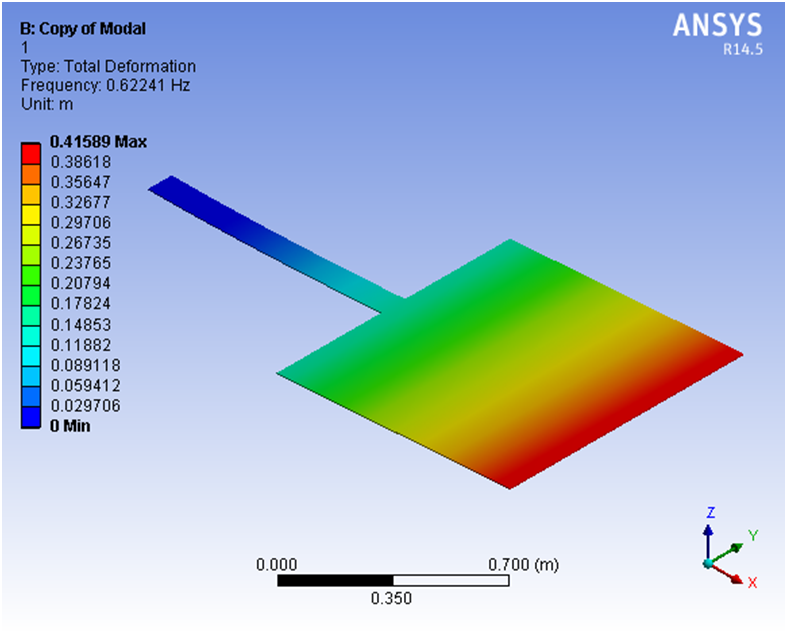} \quad
	\includegraphics[width=0.35\linewidth]{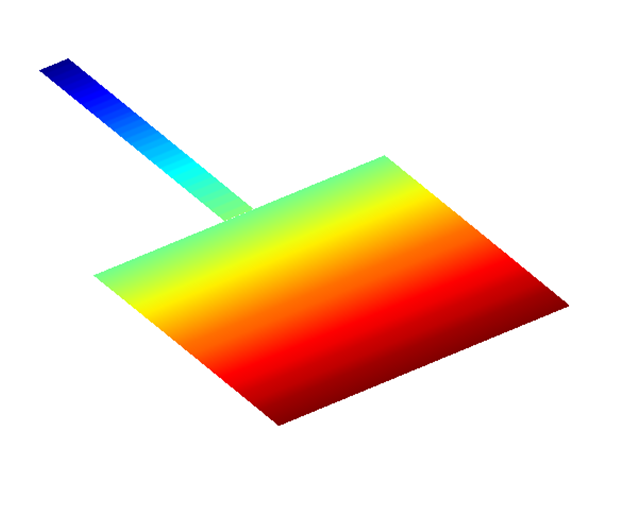}
	\thicklines
	\setlength{\unitlength}{5mm}
	\put(-10,10.5){Frequency: 0.593 Hz}
	\caption{First mode shape of transverse vibration of the beam}
	\label{fig:4}
\end{figure}

\begin{figure}[!htbp]
	\centering
	\includegraphics[width=0.5\linewidth]{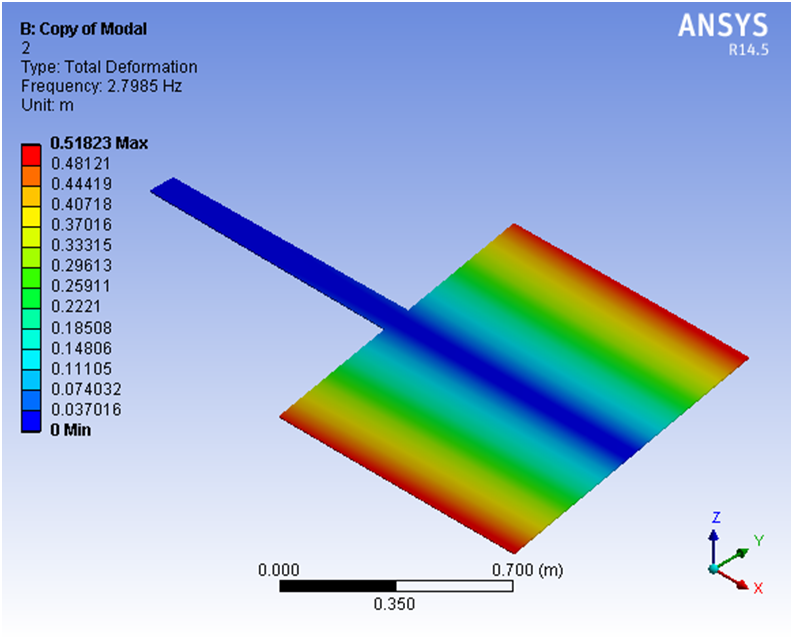} \quad
	\includegraphics[width=0.35\linewidth]{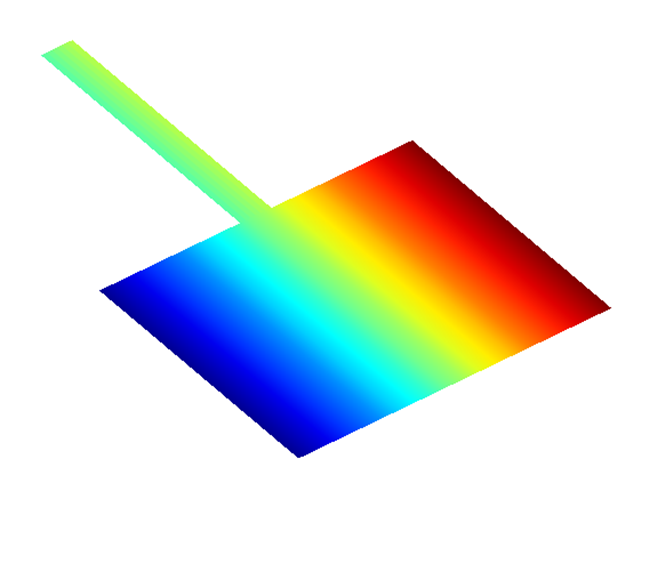}
	\thicklines
	\setlength{\unitlength}{5mm}
	\put(-10,10.5){Frequency: 2.675 Hz}
	\caption{First mode shape of torsional vibration of the beam }
	\label{fig:5}
\end{figure}

\begin{figure}[!htbp]
	\centering
	\includegraphics[width=0.5\linewidth]{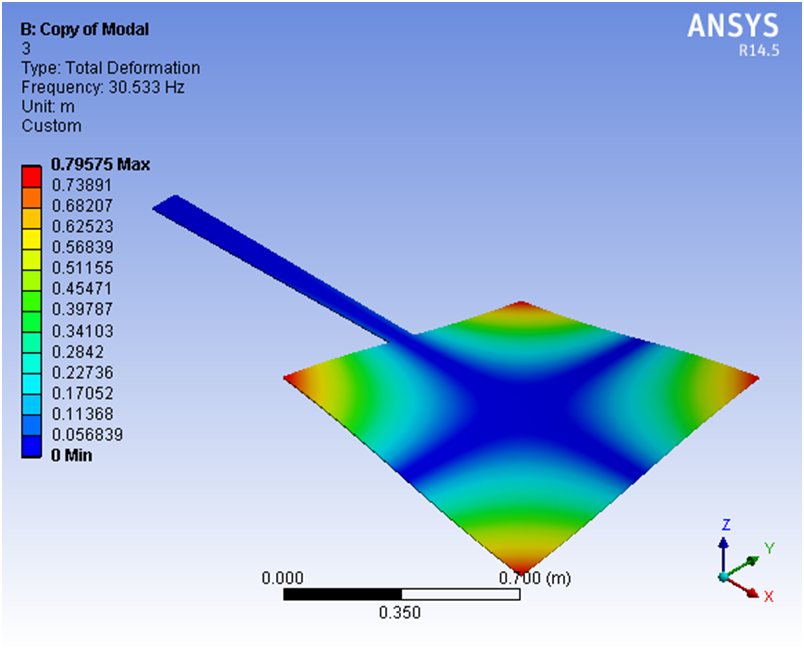} \quad
	\includegraphics[width=0.35\linewidth]{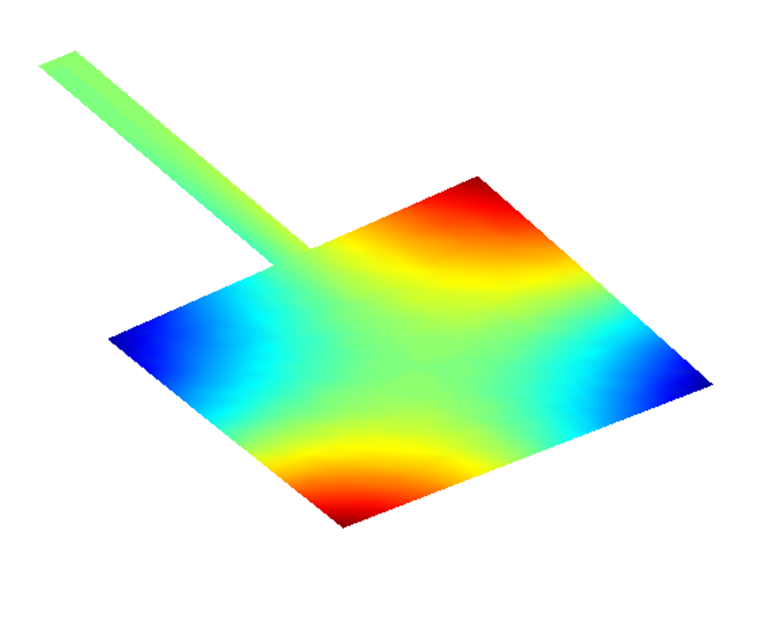}
	\thicklines
	\setlength{\unitlength}{5mm}
	\put(-10,10.5){Frequency: 32.037 Hz}
	\caption{First mode shape of torsional vibration of the beam together with first mode shape of the plate}
	\label{fig:6}
\end{figure}

\begin{figure}[!htbp]
	\centering
	\includegraphics[width=0.5\linewidth]{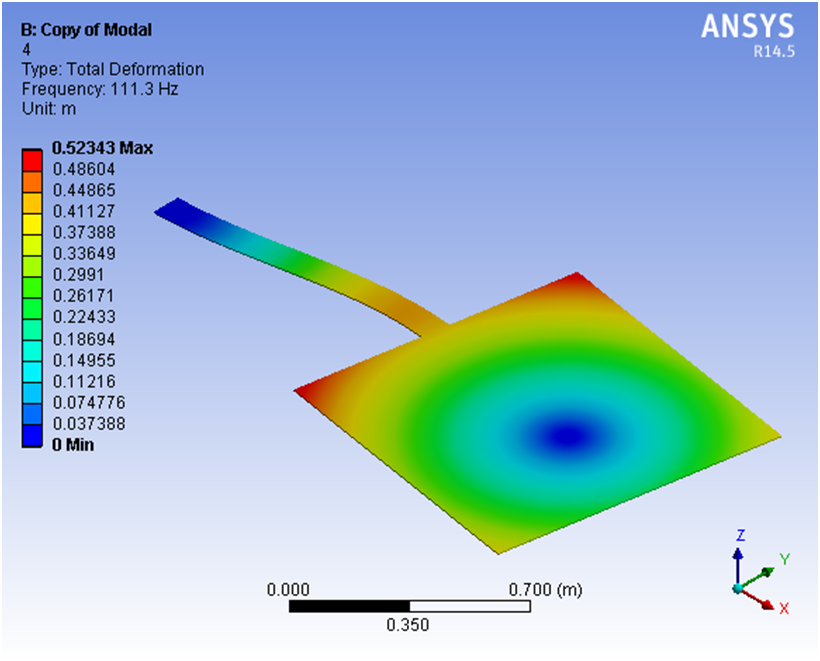} \qquad
	\includegraphics[width=0.35\linewidth]{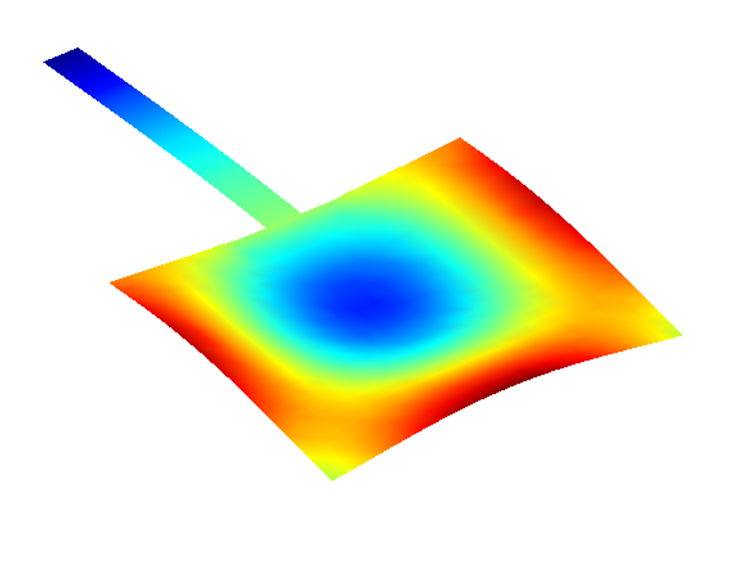}
	\thicklines
	\setlength{\unitlength}{5mm}
	\put(-10,10.5){Frequency: 118.276 Hz}
	\caption{First mode shape of transverse vibration of the beam together with second mode shape of the plate}
	\label{fig:7}
\end{figure}

\begin{figure}[!htbp]
	\centering
	\includegraphics[width=0.5\linewidth]{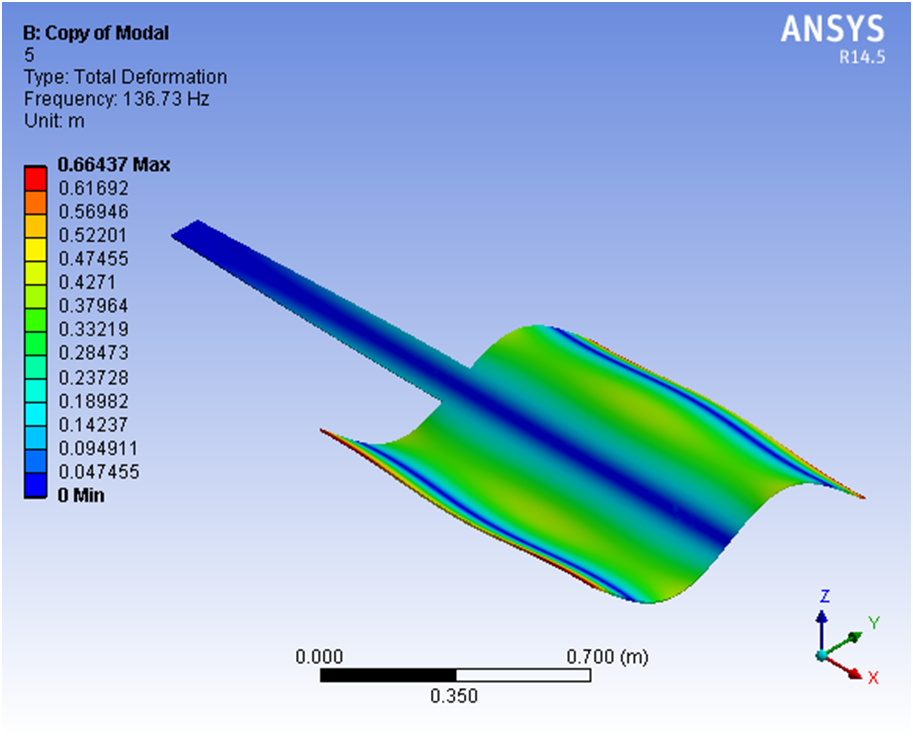} \qquad
	\includegraphics[width=0.35\linewidth]{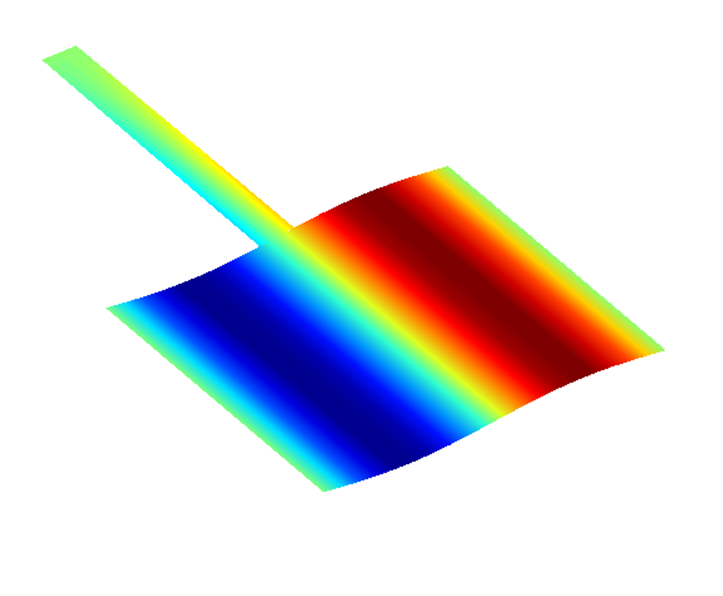}
	\thicklines
	\setlength{\unitlength}{5mm}
	\put(-10,10.5){Frequency: 147.494 Hz}
	\caption{First mode shape of torsional vibration of the beam together with forth mode shape of the plate}
	\label{fig:8}
\end{figure}

\section{Conclusion}
\label{sec:conclusion}
This study is dedicated to investigate combined torsional and transverse free vibration analysis of thin plates locally supported by an elastic beam. In particular, the performance of the GDQ method to offer solutions for practical problems with local effects is examined. Using Euler-Bernoulli assumption and classical plate theory, two coupled governing PDEs together with continuity/boundary conditions are developed based on the GDQ method. All equations are then restated as two sets of algebraic equations. The solution to the final equations leads to eigenvalues and eigenvectors that respectively indicate the system's natural frequencies and mode shapes. To assess the reliability of the proposed method, results for different well-known cases are compared by available analytical and numerical methods. This includes examples of several well-known cases of classical vibration such as cantilever beam, beam with concentrated mass, FFFF and CFFF plate. Results revealed that predictions in these case studies are in close correlation with analytical/numerical results while the accuracy can also be remarkably increased by reselecting the grid points.

Finally, predictions for the complicated system are validated by commercial FE code ANSYS, due to lack of information within the literature. Comparisons for different frequencies and mode shapes with FE results are encouraging particularly for systems with local boundary and continuity conditions. Considering the performance of the presented modeling and solution technique, it is expected to be used as a benchmark for vibration analysis of suspended plates/shells and panels in future theoretical and experimental studies.

\bibliographystyle{unsrt}  
\bibliography{Ref}

\begin{thebibliography}{10}

\bibitem{Amereh2016}
M~Amereh, MM~Aghdam, and M~Golkaram.
\newblock Design and modeling of a novel translational and angular
  micro-electromechanical accelerometer.
\newblock {\em Aerospace Science and Technology}, 50:15--24, 2016.

\bibitem{Tabak2010}
FC~Tabak, ECM Disseldorp, GH~Wortel, AJ~Katan, MBS Hesselberth, TH~Oosterkamp,
  JWM Frenken, and WM~van Spengen.
\newblock Mems-based fast scanning probe microscopes.
\newblock {\em Ultramicroscopy}, 110(6):599--604, 2010.

\bibitem{Yu2001}
Jyh-Cheng Yu and Chin-Bing Lan.
\newblock System modeling of microaccelerometer using piezoelectric thin films.
\newblock {\em Sensors and Actuators A: Physical}, 88(2):178--186, 2001.

\bibitem{Lobontiu2004}
Nicolae Lobontiu and Ephrahim Garcia.
\newblock {\em Mechanics of microelectromechanical systems}.
\newblock Springer Science \& Business Media, 2004.

\bibitem{Matsuno1996}
Fumitoshi Matsuno, Michinori Hatayama, Hideaki Senda, Tomoaki Ishibe, and
  Yoshiyuki Sakawa.
\newblock Modeling and control of a flexible solar array paddle as a
  clamped-free-free-free rectangular plate.
\newblock {\em Automatica}, 32(1):49--58, 1996.

\bibitem{Rao2007}
Singiresu~S Rao.
\newblock {\em Vibration of continuous systems}.
\newblock John Wiley \& Sons, 2007.

\bibitem{Rao1973}
G~Venkateswara Rao, IS~Raju, and TVGK Murthy.
\newblock Vibration of rectangular plates with mixed boundary conditions.
\newblock {\em Journal of Sound and Vibration}, 30(2):257--260, 1973.

\bibitem{Narita1981}
Y~Narita.
\newblock Application of a series-type method to vibration of orthotropic
  rectangular plates with mixed boundary conditions.
\newblock {\em Journal of Sound and Vibration}, 77(3):345--355, 1981.

\bibitem{Leissa1973}
Arthur~W Leissa.
\newblock The free vibration of rectangular plates.
\newblock {\em Journal of Sound and vibration}, 31(3):257--293, 1973.

\bibitem{Leissa1984}
Arthur~W Leissa and Y~Narita.
\newblock Vibrations of completely free shallow shells of rectangular planform.
\newblock {\em Journal of Sound and Vibration}, 96(2):207--218, 1984.

\bibitem{Zitnan1996}
P~{\v{Z}}it{\v{n}}an.
\newblock Vibration analysis of membranes and plates by a discrete least
  squares technique.
\newblock {\em Journal of Sound and Vibration}, 195(4):595--605, 1996.

\bibitem{Chia1985}
CY~Chia.
\newblock Non-linear vibration of anisotropic rectangular plates with
  non-uniform edge constraints.
\newblock {\em Journal of Sound and Vibration}, 101(4):539--550, 1985.

\bibitem{Leipholz1987}
HHE Leipholz.
\newblock On some dewejopments in direct methods of the caicylus of wariations.
\newblock {\em Appl Mech Rev}, 40(10):1379, 1987.

\bibitem{Donning1998}
Brian~M Donning and Wing~Kam Liu.
\newblock Meshless methods for shear-deformable beams and plates.
\newblock {\em Computer Methods in Applied Mechanics and Engineering},
  152(1-2):47--71, 1998.

\bibitem{Belytschko1994}
Ted Belytschko, Yun~Yun Lu, and Lei Gu.
\newblock Element-free galerkin methods.
\newblock {\em International journal for numerical methods in engineering},
  37(2):229--256, 1994.

\bibitem{Jang1989}
Sung~K Jang, Charles~W Bert, and Alfred~G Striz.
\newblock Application of differential quadrature to static analysis of
  structural components.
\newblock {\em International Journal for Numerical Methods in Engineering},
  28(3):561--577, 1989.

\bibitem{Cheung1989}
YK~Cheung and Wanji Chen.
\newblock Hybrid quadrilateral element based on mindlin/reissner plate theory.
\newblock {\em Computers \& structures}, 32(2):327--339, 1989.

\bibitem{Sfahani2011}
MG~Sfahani, Amin Barari, M~Omidvar, SS~Ganji, and G~Domairry.
\newblock Dynamic response of inextensible beams by improved energy balance
  method.
\newblock {\em Proceedings of the Institution of Mechanical Engineers, Part K:
  Journal of Multi-body Dynamics}, 225(1):66--73, 2011.

\bibitem{Bert1993}
Charles~W Bert, Wang Xinwei, and Alfred~G Striz.
\newblock Differential quadrature for static and free vibration analyses of
  anisotropic plates.
\newblock {\em International Journal of Solids and Structures},
  30(13):1737--1744, 1993.

\bibitem{Shu1999}
C~Shu and CM~Wang.
\newblock Treatment of mixed and nonuniform boundary conditions in gdq
  vibration analysis of rectangular plates.
\newblock {\em Engineering Structures}, 21(2):125--134, 1999.

\bibitem{Ng1999}
TY~Ng, LI~Hua, KY~Lam, and CT~Loy.
\newblock Parametric instability of conical shells by the generalized
  differential quadrature method.
\newblock {\em International Journal for Numerical Methods in Engineering},
  44(6):819--837, 1999.

\bibitem{Bellman1971}
Richard Bellman and John Casti.
\newblock Differential quadrature and long-term integration.
\newblock {\em Journal of Mathematical Analysis and Applications},
  34(2):235--238, 1971.

\bibitem{Shu1990}
C~Shu and BE~Richards.
\newblock High resolution of natural convection in a square cavity by
  generalized differential quadrature.
\newblock In {\em Proceedings of the 3rd International Conference on Advances
  in Numeric Methods in Engineering: Theory and Application, Swansea, UK},
  pages 978--985, 1990.

\bibitem{Kwak2007}
Moon~K Kwak and Sangbo Han.
\newblock Free vibration analysis of rectangular plate with a hole by means of
  independent coordinate coupling method.
\newblock {\em Journal of Sound and Vibration}, 306(1):12--30, 2007.

\bibitem{Ingber1992}
MS~Ingber, AL~Pate, and JM~Salazar.
\newblock Vibration of a clamped plate with concentrated mass and spring
  attachments.
\newblock {\em Journal of Sound and Vibration}, 153(1):143--166, 1992.

\bibitem{Ashour2004}
AS~Ashour.
\newblock Vibration of variable thickness plates with edges elastically
  restrained against translation and rotation.
\newblock {\em Thin-walled structures}, 42(1):1--24, 2004.

\bibitem{Arani2013}
Ali~Ghorbanpour Arani, Reza Kolahchi, Ali Akbar~Mosallaie Barzoki,
  Mohammad~Reza Mozdianfard, and S~Mosatafa~Noudeh Farahani.
\newblock Elastic foundation effect on nonlinear thermo-vibration of embedded
  double-layered orthotropic graphene sheets using differential quadrature
  method.
\newblock {\em Proceedings of the Institution of Mechanical Engineers, Part C:
  Journal of Mechanical Engineering Science}, 227(4):862--879, 2013.

\bibitem{Shu2001}
C~Shu, W~Chen, H~Xue, and H~Du.
\newblock Numerical study of grid distribution effect on accuracy of dq
  analysis of beams and plates by error estimation of derivative approximation.
\newblock {\em International Journal for Numerical Methods in Engineering},
  51(2):159--179, 2001.

\bibitem{Shu1997}
C~Shu and H~Du.
\newblock A generalized approach for implementing general boundary conditions
  in the gdq free vibration analysis of plates.
\newblock {\em International Journal of Solids and Structures}, 34(7):837--846,
  1997.

\end{thebibliography}

\end{document}